\documentclass[11pt]{amsart}
 
\usepackage{amscd}
 
\newcommand{\ot}{\otimes}

\renewcommand{\P}{\mathbb{P}}

\newcommand{\Z}{{\mathbb{Z}}}        

\newcommand{\QZt}{\mathbb{Q}/\mathbb{Z}(2)}

\newcommand{\Zm}[1]{\Z/{#1}\Z}

\newcommand{\La}{\Lambda}
\newcommand{\la}{\lambda}
\newcommand{\D}{\Delta}

\newcommand{\G}{{\Gamma}}       

\newcommand{\oddots}{{\mathinner{\mkern1mu\raise1pt\vbox{\kern7pt\hbox{.}}\mkern2mu\raise4pt\hbox{.}\mkern2mu\raise7pt\hbox{.}\mkern1mu}}}

\newcommand{\basemu}{\boldsymbol{\mu}}
\newcommand{\mmu}[1]{\basemu_#1}     

\DeclareMathOperator{\Spin}{Spin}           

\newcommand{\Gm}{\mathbb{G}_m}

\newcommand{\oE}{\ensuremath{^1\!E_6}}
\newcommand{\dE}{\ensuremath{^2\!E_6}}

\DeclareMathOperator{\Spec}{Spec}

\DeclareMathOperator{\Gal}{Gal}  

\DeclareMathOperator{\im}{im}

\DeclareMathOperator{\chr}{char}

\DeclareMathOperator{\cores}{cor}

\DeclareMathOperator{\End}{End}

\DeclareMathOperator{\Aut}{Aut}

\DeclareMathOperator{\Iso}{Iso}

\newcommand{\iso}{\xrightarrow{\sim}}

\newcommand{\ra}{\rightarrow}

\newcommand{\qf}[1]{{\langle{#1}\rangle}} 
\newcommand{\pff}[1]{{\langle\!\langle{#1}\rangle\!\rangle}} 
\newcommand{\qpf}[1]{{\langle\!\langle{#1}]]}} 
 
\newcommand{\kx}{k^{\times}}
\newcommand{\kxii}{k^{\times} / k^{\times 2}}

\renewcommand{\sp}{{\mathrm{sp}}}

\renewcommand{\P}{\mathbb{P}}
 
\newcommand{\an}{{\mathrm{an}}}
 
\newcommand{\lbar}{\overline{\ell}}
 
\newcommand{\oEiso}{\ensuremath{^1\!E^{28}_{6,2}}}    

\newcommand{\darkradE}{0.15}
\newcommand{\lradE}{0.4}
 
\usepackage{amsthm}

\newtheorem{thm}[equation]{Theorem}
\newtheorem{lem}[equation]{Lemma}
\newtheorem{cor}[equation]{Corollary}
\newtheorem{prop}[equation]{Proposition}

\newtheorem*{thm*}{Theorem}
\newtheorem*{prop*}{Proposition}
\newtheorem*{cor*}{Corollary}
\newtheorem*{lem*}{Lemma}

\theoremstyle{definition} %

\newtheorem*{defn*}{Definition}

\theoremstyle{remark} %
\newtheorem{rmk}[equation]{Remark}

\newtheorem*{rmk*}{Remark}

\numberwithin{equation}{section}

\makeatletter
{\smallskip \refstepcounter{equation}\noindent{\textbf \theequation. }{{\textbf{#1.}}}}%
{\smallskip \global\@ignoretrue}
\makeatother

\makeatletter
{\smallskip \refstepcounter{equation}{\sc \theequation}{\sc (#1).}}%
{\smallskip \global\@ignoretrue}
\makeatother

\makeatletter
{\smallskip \refstepcounter{equation}{\sc \theequation.}{\sl{ #1.}}}%
{\smallskip \global\@ignoretrue}
\makeatother

\makeatletter
\newenvironment{borel*}%
{\smallskip \refstepcounter{equation}\noindent{\textbf \theequation.}}%
{\global\@ignoretrue}
\makeatother

\begin{document}

\title{Totaro's Question for $G_2$, $F_4$, and $E_6$}
\author{Skip Garibaldi}
\address{Department of Mathematics \& Computer Science, Emory University, 
Atlanta, Georgia 30322, USA}
\email{skip@member.ams.org}
\urladdr{http://www.mathcs.emory.edu/{\textasciitilde}skip/}
 
\author{Detlev Hoffmann}
\address{Division of Pure Mathematics, School of Mathematical Sciences, 
University of Nottingham, University Park, Nottingham NG7 2RD, UK}
\email{Detlev.Hoffmann@Nottingham.ac.uk}              

\begin{abstract}
In a 2004 paper, Totaro asked whether a $G$-torsor $X$ that has a 
zero-cycle of degree $d>0$ will necessarily have a closed \'etale point
of degree dividing $d$, where $G$ is a connected algebraic group.  
This question is closely related to several conjectures regarding
exceptional algebraic groups. Totaro gave a positive answer to his 
question in the following cases: $G$ simple, split, and of type $G_2$, 
type $F_4$, or simply connected of type $E_6$. We extend the list of 
cases where the answer is ``yes" to all groups of type $G_2$ and some 
nonsplit groups  of type $F_4$ and $E_6$. No assumption on the characteristic 
of the base field is made. The key tool is a lemma regarding linkage of Pfister forms.
\end{abstract}

 
\maketitle
 
\setcounter{equation}{1}

\nocite{Baeza2}

For certain linear algebraic groups $G$ over a field $k$ and certain 
homogeneous $G$-varieties $X$, Totaro asked in \cite{Tot:E8}:
\begin{equation} \label{ques}
\parbox{4in}{\emph{If $X$ has a zero-cycle of degree $d > 0$, does $X$ 
necessarily have a closed \'etale
point of degree dividing $d$?}}
\end{equation}
(This question is closely related to earlier questions raised by Veisfeiler, 
Serre, and Colliot-Th\'el\`ene.)  This question can be rephrased as:
\begin{equation}
\tag{\ref{ques}$'$}
\parbox{4in}{\emph{If $X$ has a (closed) point over finite extensions $K_1, K_2, \ldots, K_n$ of $k$, does $X$ necessarily have a point over a separable extension of $k$ of degree dividing every $[K_i:k]$?}}
\end{equation}
The purpose of this note is to prove the following theorem.

\begin{thm} \label{MT}
The answer to \eqref{ques} is ``yes'' when $X$ is a $G$-torsor and $G$ is
   \begin{itemize}
   \item of type $G_2$,
   \item reduced of type $F_4$, or
   \item simply connected of type $^1\!E^0_{6,6}$ (split) or of type \oEiso.
   \end{itemize}
\end{thm}

A \emph{$G$-torsor} is a principal $G$-bundle over $\Spec k$.   
Torsors are often called principal homogeneous spaces, as in \cite{SeCG}. 
It is possible that the answer to \eqref{ques} is ``yes'' whenever $G$ is semisimple and $X$ is a $G$-torsor (and, in particular, is affine); no counterexamples are known.  In contrast, for $X$ projective, there are examples where the answer is ``no''  even when $d = 1$, see \cite{Flo} and \cite{Par}.

Every group of type $F_4$ is the group of automorphisms of some uniquely 
determined Albert $k$-algebra
$J$.  We say that the group is \emph{reduced} if the algebra $J$ is reduced, 
i.e., if $J$ is not a division algebra.  (Recall that Albert algebras are 27-dimensional 
Jordan algebras.  For a survey of the state of the art in 1992, see
\cite{PR}.  Many useful and characteristic-free results  about Albert 
algebras may be found in \cite{Ptr:struct}.)

The notation ``\oEiso'' is from Tits' classification \cite[p.~58]{Ti:Cl}.  
Such groups have Tits index
\setlength{\unitlength}{.75cm}
\begin{equation}  \label{e6.ind} 
\begin{picture}(5,2)
    \multiput(0.5,0.5)(1,0){5}{\circle*{\darkradE}}
    \put(2.5,1.5){\circle*{\darkradE}}

    \put(0.5,0.5){\line(1,0){4}}
    \put(2.5,1.5){\line(0,-1){1}}

    \put(0.5,0.5){\circle{\lradE}}
    \put(4.5,0.5){\circle{\lradE}}
\end{picture}
\end{equation}

\smallskip

Totaro proved Theorem \ref{MT} under the additional assumption that $G$ is split, 
see \cite[5.1]{Tot:E8}.  The main tool that allows us to extend his result is 
our Proposition \ref{quad.prop} and its analogue for cohomology, Cor.~\ref{coho}; 
they assert the existence of field extensions of small dimension that kill a 
difference of symbols in $H^p(k, \Zm2(p-1))$.  The quadratic form theory that 
underlies Cor.~\ref{coho} occupies half of this paper
and concerns mainly the question of so-called linkage of Pfister forms
(in all characteristics, including $2$).

The last section gives some comments on extending Theorem \ref{MT} to include 
other exceptional groups.
Except for this final section, we make no assumptions on the base field $k$.

\subsection*{Notation} Throughout, we write $\mmu{n}$ for the 
group scheme of $n$-th roots of unity, where $n$ is a natural number.  
We write $H^p(k, \Zm{n}(p - 1))$ for the abelian group defined in \cite[App.~A]{MG};
when $n$ is not divisible by the characteristic of $k$, 
it is the Galois cohomology group $H^p(k, \mmu{n}^{\otimes (p - 1)})$.

\section{Quadratic form preliminaries}

For readers who are not very familiar with the algebraic theory
of quadratic forms we will first
present a few definitions and facts.  
Other basic results and terminologies that
we use without further reference can all be found in Scharlau's book
\cite{Sch} for characteristic $\neq 2$, and in Baeza's book \cite{Baeza}
and the article \cite{HL} by Laghribi and the second author for the
characteristic $2$ case.
In particular, we will freely use certain properties of
Pfister forms, Pfister neighbors, and function fields of quadratic forms.

If $\alpha$ and $\beta$ are quadratic forms,
then $\beta$ is called a \emph{subform} of $\alpha$ 
(denoted $\beta\subset\alpha$), if there exists a form
$\gamma$ such that $\alpha\cong\beta\perp\gamma$, and 
$\alpha$ is said to \emph{dominate} $\beta$ (denoted
$\beta\prec\alpha$) if $\beta$ is isometric to the
restriction of $\alpha$ to some subspace of the underlying vector space
of $\alpha$. 

The Witt index $i_W(\phi)$ of a  quadratic form $\phi$ is the
maximal number $m$ such that $\phi$ contains an orthogonal sum of
$m$ hyperbolic planes as a subform. Furthermore,
if $\beta$ is a bilinear form (in characteristic $2$), then we denote by
$\beta_q$ the totally singular quadratic form defined by $\beta_q(x)=\beta
(x,x)$.

Every nonsingular quadratic form $\phi$ over $k$ has attached to it
an invariant $e_1(\phi)$ living in the Galois cohomology group $H^1(k,\Zm2)$.
In characteristic not $2$, we have $H^1(k,\Zm2)=\kxii$ and $e_1(\phi)=d_{\pm}(\phi)$ 
is the signed discriminant.  If $\chr (k)=2$, then $H^1(k,\Zm2)=k/\wp (k)$
and $e_1(\phi )=\Delta (\phi)$ is the Arf-invariant.
Note that in all characteristics $H^1(k,\Zm2)$ is nothing else but the group 
$H^1(k,\Zm2(0))$ from above.

Recall that in characteristic $\neq 2$ (resp.\ characteristic $2$),
a product of $p$ binary
quadratic (resp.\ bilinear) forms $\qf{1,-a_i}$ is called a
$p$-fold Pfister form (resp.\ $p$-fold bilinear Pfister form):
$\pff{a_1,\cdots,a_p}=\qf{1,-a_1}\otimes\cdots\otimes\qf{1,-a_p}$.
These forms generate
additively the ideal $I^pk$, the $p$-th
power of the fundamental ideal of even-dimensional forms in $Wk$, the Witt
ring of quadratic forms (resp.\ bilinear forms) over $k$.

In characteristic $2$, a $p$-fold quadratic Pfister form is a product
of a $(p-1)$-fold bilinear Pfister form $\pff{a_1,\cdots ,a_{p-1}}$
and a nonsingular quadratic form $[1,a_p]=X^2+XY+a_pY^2$, $a_p\in k$:
$\qpf{a_1,\cdots,a_{p-1},a_p}= \pff{a_1,\cdots ,a_{p-1}}\otimes [1,a_p]$.
The $p$-fold quadratic Pfister forms generate the $Wk$-submodule
$I^p_qk =(I^{p-1}k)W_qk$ of the Witt group $W_qk$ of nonsingular
quadratic forms over $k$.

Here, a {\em $p$-Pfister form} will always mean either a
$p$-fold Pfister form in characteristic $\neq 2$
or a $p$-fold quadratic Pfister form if the characteristic is $2$, and we
will only consider the case $p\geq 1$.
Also, we will use in characteristic $\neq 2$ the notation $I^p_qk$ instead
of $I^pk$ for uniformity's sake, and we will denote the set of
isometry classes of $p$-Pfister forms over $k$ by $P_pk$.

A quadratic form $\phi$ is a {\em Pfister neighbor} of a Pfister form $\pi$ if
$a\phi\subset\pi$ ($\chr(k)\neq 2$) resp.\ $a\phi\prec\pi$ ($\chr(k)=2$) 
for some $a\in\kx$, and $2\dim\phi >\dim\pi$.  

If $\phi$ is a Pfister neighbor of $\pi$ and $E/k$ is any field
extension then $i_W(\phi_E)>0$ if and only if $\pi_E$ is isotropic if and only if
$\pi_E$ is hyperbolic.  
If any of these equivalent conditions hold, then
$E(\phi)/E$ and $E(\pi)/E$ are purely transcendental, and any
form $\psi$ over $k$ which is isotropic over $E(\phi)$ or $E(\pi)$ is so
already over $E$. Applied to $E=k(\phi )$ resp.\ $E=k(\pi)$, this shows
that a form $\psi$ over $k$ is isotropic over $k(\pi)$ if and only if it is
isotropic over $k(\phi)$.

In characteristic $\neq 2$, it is well-known that if $\phi$ and $\pi$ are
Pfister forms with $\pi\subset\phi$,
then there exists a Pfister form $\rho$ such that
$\phi$ is isometric to $\rho\otimes\pi$.
The analogous result in characteristic $2$ is the following.

\begin{lem} \label{B} Suppose $\chr(k)=2$. Let $\phi$ be a Pfister form
over $k$ and let $x$ be in $\kx$.
\begin{enumerate}
\item[(i)] Let $\pi$ be a Pfister form with $x\pi\subset\phi$.
Then there exists a bilinear Pfister form $\rho$ such that
$\phi\cong\rho\otimes\pi$.
\item[(ii)]  Let $\rho$ be a bilinear Pfister form such that
$x\rho_q\prec\phi$.  Then
there exists a Pfister form $\pi$ such that $\phi\cong\rho\otimes\pi$.
\end{enumerate}
\end{lem}

\begin{proof}  We assume that $\phi$ is anisotropic and leave
the (much simpler) case of $\phi$ being isotropic (and hence hyperbolic)
to the reader.  We may also assume that $x=1$ because of the roundness
of Pfister forms. (Recall that a form is called \emph{round} if the set of 
nonzero elements represented by the form is exactly the group of similarity
factors of the form.)  Indeed, since $\pi$ resp.\ $\rho_q$ represent $1$,
the hypotheses in (i) resp.\ (ii) imply that $\phi$ represents $x$,
hence, the roundness of $\phi$ implies $\phi\cong x\phi$ and thus 
$\pi\subset\phi$ resp. $\rho_q\prec\phi$.

(i) Let $\mu$ be a quadratic form such that $\phi\cong\mu\perp\pi$.
Now $\pi =\qf{1}\otimes\pi$ and $\qf{1}$ is a
$0$-fold bilinear Pfister form.  So let $r\geq 0$ be maximal such that
there exists an $r$-fold bilinear Pfister form $\rho$ and a quadratic
form $\eta$ with $\phi\cong\eta\perp (\rho\otimes\pi)$.

Suppose $\dim\eta>0$, and let $a\in\kx$ be represented by $\eta$.
Then $\gamma'=\qf{a}_q\perp\rho\otimes\pi$ is dominated by $\phi$.
In particular, $\phi$ becomes isotropic and hence hyperbolic
over the function field $k(\gamma')$.  But $\gamma'$ is a Pfister
neighbor of the Pfister form $\gamma=\pff{a}\otimes\rho\otimes\pi$,
hence $\phi$ becomes isotropic and hence hyperbolic over $k(\gamma)$.
Since $\gamma$ and $\phi$ represent $1$, it follows that there exists
a form $\eta'$ such that $\phi\cong\eta'\perp\gamma$ (see, e.g.,
\cite[Th.~4.2]{HL}).
We get a contradiction to the maximality of $r$ as
$\gamma\cong\rho'\otimes\pi$ with
$\rho'=\pff{a}\otimes\rho$ a bilinear Pfister form of fold $r+1$.
Hence $\dim\eta =0$ and $\phi\cong\rho\otimes\pi$.

(ii) Write $\rho_q\cong\qf{1}_q\perp\rho'_q$ (with $\rho'_q$ totally
singular).  Since $\rho_q$ is dominated by $\phi$, there exists $b\in k$
such that $\gamma'=[1,b]\perp\rho'_q$ is also dominated by $\phi$
(see \cite[Cor.~3.3, Lemma 3.5]{HL}).  Therefore, $\phi$ becomes isotropic and hence
hyperbolic over $k(\gamma')$.  But $\gamma'$ is a Pfister neighbor of
the Pfister form $\gamma=\rho\otimes [1,b]$.  It follows that $\phi$ becomes
hyperbolic over $k(\gamma)$, and as in the proof of (i), we find
a bilinear Pfister form $\nu$ such that $\phi\cong\nu\otimes\rho\otimes [1,b]$.
The claim follows by putting $\pi =\nu\otimes [1,b]$.
\end{proof}

\begin{rmk}  Part (i) in the above lemma is nothing else but
\cite[Th.~4.4]{Baeza}.  Our proof is quite different from Baeza's
original one in that it uses the machinery of function
fields of quadratic forms.  Part (ii) is basically the analogous result
with the roles of bilinear and quadratic Pfister forms reversed.
\end{rmk}

\begin{lem} \label{EL}  Let $\phi$ and $\psi$ be anisotropic Pfister forms of
fold $p$ and $q$ respectively.  Then there exists an integer $m\geq 0$
such that $i_W(\phi\perp -\psi)=2^m$.

Furthermore, if $m>0$, then $m$ is maximal for the property such that there exist
an $m$-Pfister form $\pi$ and Pfister forms (resp.\ bilinear Pfister forms
if $\chr(k)=2$) $\rho$ and $\sigma$ of fold $p-m$ and $q-m$, respectively,
with $\phi\cong\rho\otimes\pi$ and $\psi\cong\sigma\otimes\pi$.
\end{lem}

\begin{proof}  If $\chr(k)\neq 2$, this is essentially \cite[Prop.4.4]{El:pf}.
In characteristic $2$, the result is due to Faivre \cite{Faivre}, but
we include a proof of this case for the reader's convenience.

So let us assume that $\chr(k)=2$.
Clearly, $\phi$ and $\psi$ represent $1$, so $i_W(\phi\perp -\psi)\geq 1$.
Also, if $\pi$ is an $m$-Pfister dividing both $\phi$ and $\psi$, then
the roundness of Pfister forms implies that $\pi$ is a subform of
both $\phi$ and $\psi$, therefore 
$$i_W(\phi\perp-\psi)\geq i_W(\pi\perp -\pi)=
\dim\pi = 2^m.$$ 
It thus suffices to consider the case $i_W(\phi\perp -\psi)\geq 2$.
We have to show that the Witt index is equal to the maximal
dimension of a Pfister form $\pi$ which is a subform
of both $\phi$ and $\psi$.
The existence of the bilinear Pfister forms $\rho$ and $\sigma$ then
follows from Lemma \ref{B}(i).

We first show that
there exists a common $1$-Pfister form as subform of both $\phi$ and $\psi$.

Now $i_W(\phi\perp -\psi)\geq 2$ implies that there is a
binary form $\beta$ (possibly singular) which is dominated
both by $\phi$ and $\psi$ (see \cite[Cor.~3.13]{HL}), and we may assume
that $\beta$ represents $1$ as $\phi$ and $\psi$ represent $1$.
If $\beta$ is isometric to $[1,b]$ for some $b\in k$, then we are done.
Otherwise, $\beta$ is isometric to $\qf{1,b}_q$ for some $b\in \kx$.  
In this case, we necessarily
have $p,q\geq 2$ for dimension reasons and since $\phi$ and $\psi$ are
nonsingular.  We then apply Lemma \ref{B}(ii) to find bilinear Pfister forms
$\rho'$ and $\sigma'$, and  $c,d\in k$ such that
$\phi\cong\pff{b}\otimes\rho'\otimes [1,c]$ and
$\psi\cong\pff{b}\otimes\sigma'\otimes [1,d]$.  But then there exist
$e,f,g\in k$ such that $\qpf{b,c}\cong\qpf{e,g}$ and
$\qpf{b,d}\cong\qpf{f,g}$ (see, e.g., \cite{Lam:link} for the equivalent
``linkage'' result for quaternion algebras).
Clearly, $[1,g]$ is the desired form.

Having found a $1$-Pfister form as subform of both $\phi$ and $\psi$,
let now $\pi$ be an $m$-Pfister form which is a subform of both
$\phi$ and $\psi$ and with $m$ maximal.  Thus, we can write
$\phi\cong\phi'\perp\pi$ and $\psi\cong\psi'\perp\pi$,
and we have
$$i_W(\phi\perp -\psi)=i_W(\phi'\perp -\psi')+2^m.$$
If $i_W(\phi\perp -\psi)>2^m$ then $\phi'\perp -\psi'$ is isotropic,
hence $\phi'$ and $\psi'$ represent a common element $a\in \kx$, and 
an argument
as in the proof of Lemma \ref{B}(i) then shows that the $(m+1)$-Pfister
$\pff{a}\otimes\pi$ is a subform of both $\phi$ and $\psi$,
a contradiction to the maximality of $m$.
This completes the proof.
\end{proof}

In the above lemma, we restricted ourselves to the
case of anisotropic Pfister forms.  It is obvious that if, say, $\phi$
is isotropic and hence hyperbolic, then $i_W(\phi\perp -\psi)=2^{p-1}$
or $2^{p-1}+2^{q-1}$, depending on whether $\psi$ is anisotropic or not.
In this case, $\phi$ divides $\psi$ if $\psi$ is hyperbolic and $q\geq p$,
and $\psi$ divides $\phi$ if $p>q$ or $p=q$ and $\psi$ is hyperbolic.
If $\psi$ is anisotropic and $q\geq p$, then any $(p-1)$-Pfister dividing
$\psi$ will also divide the hyperbolic form $\phi$, but there is
obviously no $p$-Pfister dividing both $\phi$ and $\psi$.

\begin{lem} \label{DHlem}
Let $\phi$, $\psi$ be $p$-Pfister forms over $k$, with $p \ge 2$.
Suppose that there is an odd-degree extension $L$ of $k$ such that
$\phi$ and $\psi$ are isometric over a quadratic extension of $L$.
Then there exists a $1$-Pfister $\alpha$ dividing both $\phi$ and $\psi$.
In particular, $\phi$ and $\psi$ are both hyperbolic over $k$ or some
separable quadratic extension of $k$.
\end{lem}

If $p = 1$, let $c=e_1(\phi\perp -\psi)\in H^1(k,\Zm2)$.
We then have that $\phi$ and $\psi$ are isometric if and only if
$c\in H^1(k,\Zm2)$ is trivial.
Hence $\phi$ and $\psi$ are isometric over $k(\sqrt{c})$ ($\chr (k)\neq 2$)
resp.\ $k(\wp^{-1}(c))$ ($\chr (k)=2$).
But $\phi$ and $\psi$ need not be divisible by the same
$1$-Pfister.

The statement of the lemma might also be expressed by saying that
the two Pfister forms have a common right slot.  In other words, 
there exist $c\in k$
and $a_1,\ldots ,a_{p-1},b_1,\ldots ,b_{p-1}\in \kx$ such that the following
holds: 

If $\chr(k)\neq 2$, we have
$\phi\cong\pff{a_1,\ldots ,a_{p-1},c}$, $\psi\cong\pff{b_1,\ldots,b_{p-1},c}$,
and the extension is $k(\sqrt{c})/k$.

If $\chr(k)=2$, we have
$\phi\cong\qpf{a_1,\ldots ,a_{p-1},c}$, $\psi\cong\qpf{b_1,\ldots,b_{p-1},c}$,
and the extension is $k(\wp^{-1}(c))/k$.

Of course, the order of the slots only matters in characteristic $2$.

\begin{proof}
Suppose first that $\phi$ is isotropic and hence hyperbolic, and let
$b_1,\ldots $, $b_{p-1}\in \kx$, $c\in k$ be such that
$\psi\cong\pff{b_1,\ldots,b_{p-1},c}$ ($\chr (k)\neq 2$) resp.
$\psi\cong\qpf{b_1,\ldots,b_{p-1},c}$ ($\chr(k)=2$).  Since $\phi$
is hyperbolic, we have
$\phi\cong\pff{-1,\ldots ,-1,c}$ resp.\ $\phi\cong\qpf{-1,\ldots ,-1,c}$,
and the result is obvious.

So let us assume that both $\phi$ and $\psi$ are anisotropic.
Set $q$ to be the anisotropic part of $\phi\perp - \psi$. Let $m$ be maximal
such that $\phi$ and $\psi$ contain a common $m$-Pfister as subform.
By Lemma \ref{EL}, we have
   \[
   \dim q = 2^{p+1}-2i_W(\phi\perp - \psi)=2^{p+1} - 2^{m + 1},
   \]
and it suffices to show that $m$ is not zero; clearly we may assume that 
$m < p$, i.e., that $\dim q> 0$.

By Springer's Theorem for odd-degree extensions
(see, e.g., \cite[2.5.3]{Sch}), $q$ is $L$-anisotropic.
Now $q$ is hyperbolic over a quadratic extension $E$ of $L$.  We consider
three cases:
\begin{enumerate}
\item[(i)] $\chr(k)=2$ and $E=L(\sqrt{b})$ for some $b\in L^\times$
($E/L$ is inseparable);
\item[(ii)] $\chr(k)\neq 2$ and $E=L(\sqrt{b})$ for some $b\in L^\times$;
\item[(iii)] $\chr(k)=2$ and $E=L(\wp^{-1}(b))$ for some $b\in L^\times$
($E/L$ is separable);
\end{enumerate}
In cases (i) and (ii), there exists a nonsingular quadratic form 
$\eta$ over $L$ 
such that   
$$q_L\cong\pff{b}\otimes\eta$$
(see \cite[2.5.2]{Sch} in case (i) and \cite{Ahmad} in case (ii)).
In case (iii), there exists a bilinear form $\eta$ over $L$ such that
$$q_L\cong\eta\otimes [1,b]$$
(see \cite{Ba:Teil}). 

In case (i), the nonsingularity of $\eta$ implies that $\dim\eta \equiv
0\bmod 2$, hence $\dim q\equiv 0\bmod 4$ and thus $m\geq 1$.

Suppose we are in case (ii) or (iii), and put 
$\beta =\pff{b}$ or $\beta =[1,b]$, respectively. Note that
$\beta$ is obviously anisotropic, and that 
$q_L\in I^2_qL$ as $\phi$ and $\psi$ are $p$-Pfisters with $p\geq 2$.
If $\dim\eta$ is odd, then there exists $c\in L^\times$ such that
$$0\equiv q_L\equiv c\beta\bmod I^2_qL\ ,$$
but this contradicts the Arason-Pfister Hauptsatz which states that
an\-iso\-tro\-pic forms in $I^n_qL$ are of dimension $\geq 2^n$.
Hence, again $\dim\eta\equiv0\bmod 2$ and therefore $m\geq 1$.
\end{proof}

The following field-theoretic lemma is folklore but we include a proof
for the reader's convenience.

\begin{lem} \label{field.lem}
Let $K/k$ be a finite separable extension of degree $p^nm$ where $p$
is a prime not dividing $m$.  Then there exist finite separable extensions
$E/L/k$ such that $K\subset E$, $[E:L]=p^n$,  and $p$ does not
divide $[L:k]$.
\end{lem}
\begin{proof}  Let $N/k$ be any finite Galois extension containing $K$, and
let $S'$ be any $p$-Sylow of $\Gal (N/K)$.  Since $\Gal (N/K)$ is a subgroup
of $\Gal (N/k)$,
there is a $p$-Sylow $S$ of $\Gal (N/k)$ containing $S'$.  Let $L=N^S$ and
$E=N^{S'}$.  Clearly, $E$ contains $K$, $p$ does not divide $[L:k]$,
and by comparing degrees we have $[E:L]=p^n$.
\end{proof}

\begin{prop} \label{quad.prop}
Let $\phi$ and $\psi$ be $p$-Pfister forms over a field $k$.  
Suppose that there exist finite extensions $K_1,
\ldots, K_n$ of $k$ such that $\phi$ and $\psi$ are isometric over 
$K_i$ for all $i$.  Then there is a Galois extension $K$ of $k$ such that $
\phi$ and $\psi$ are $K$-isometric and $[K:k]$ divides
   \[
   g:=\gcd\{4, [K_1:k], [K_2:k], \ldots, [K_n: k]\}.
   \]
If $p \ge 2$ and $g$ is even, then $K$ may be chosen to kill both $\phi$
and $\psi$.  If $g=1$, then there exists a Galois extension of degree dividing $2$
that kills both $\phi$ and $\psi$.
\end{prop}

\begin{proof}
If $\phi =\pff{\cdots ,a}$, $\psi =\pff{\cdots ,b}$ ($\chr (k)\neq 2$) resp.\
$\phi =\qpf{\cdots ,a}$, $\psi =\qpf{\cdots ,b}$ ($\chr(k)=2$), then
$k(\sqrt{a},\sqrt{b})$ resp.\ $k(\wp^{-1}(a),\wp^{-1}(b)$) is a Galois
extension of degree dividing $4$ which obviously kills $\phi$ and
$\psi$ simultaneously.

So we may assume that the gcd $g$ in the statement is $1$ or $2$.
If $g=1$, then $[K_i : k]$ 
is odd for some $i$, hence $\phi$ and $\psi$ are $k$-isomorphic by Springer's
theorem.  In this case, we may choose $a=b$ above and thus get the
quadratic Galois extension which kills both forms.

Finally, suppose that $g=2$.  In this case, we have $[K_i : k] \equiv 2 
\bmod 4$ for some $i$ which we fix.

Let us first consider the case $\chr(k)\neq 2$. We then can write $K_i$ as a
separable  extension $K_i'$ followed by a purely inseparable one which
necessarily is of odd degree. Then $[K_i' : k] \equiv 2 
\bmod 4$, and using Springer's theorem,  we may then assume that $K_i$ is
separable over $k$ by replacing $K_i$ by $K_i'$
if necessary.   By Lemma \ref{field.lem}, there
are separable extensions $E/L/k$ such that $K_i\subset E$, 
$[E:L]=2$, and $[L:k]$ is odd.  Since $K_i\subset E$,
we have that $\phi$ and $\psi$ are isometric over $E$, and we
conclude by applying Lemma \ref{DHlem}.

Now suppose $\chr(k)=2$.  If $K_i/k$ is separable, the same
reasoning as above applies.  If it is inseparable, then this is
only possible for reasons of degree and characteristic if there is an
intermediate field $L$ such that $L/k$ is of odd degree
(and thus separable) and $K_i/L$ is inseparable of degree $2$.
Again, Lemma \ref{DHlem} yields the desired conclusion.
\end{proof}

The above proposition allows an interpretation in terms of 
symbols in the group $H^p(k, \Zm2 (p-1))$
mentioned in the introduction.  Recall that in characteristic $\neq 2$,
$H^p(k, \Zm2 (p-1))$ is just the usual Galois cohomology group
with $\bmod\, 2$ coefficients, $H^p(k,\Zm2 )$, and that there is a well-defined
injective map $P_pk\to H^p(k,\Zm2 )$ sending $\pff{a_1,\cdots ,a_p}$ to
the symbol  $(a_1)\cup\cdots\cup (a_p)$.
This map is well-defined, see Elman-Lam \cite[3.2]{El:pf} or
Arason \cite[1.6]{Ar:th}, and it is injective
as a direct consequence of Voevodsky's proof of the Milnor conjecture.

As already mentioned, in characteristic $2$, we have $H^1(k, \Zm2 (0))=k/\wp (k)$.
One can then show that $H^p(k, \Zm2 (p-1))$ is nothing else but 
the group $h_p$ defined in \cite{ABa}, namely the quotient of the group
$(K_{p-1}k/2K_{p-1}k)\otimes H^1(k, \Zm2 (0))$ by the subgroup $\mathcal{R}_p$
generated by elements of type 
$$\{a_1,\cdots ,a_{p-1}\}\otimes (b)\in 
(K_{p-1}k/2K_{p-1}k)\otimes H^1(k, \Zm2 (0))$$
such that $a_i$ is a norm of the extension $k(\wp^{-1}(b))/k$ for some $i$
(or, equivalently, that $a_i$ is represented by the quadratic form $[1,b]$).
Again, we have a well-defined and injective map $P_pk\to h_p$ sending
$\qpf{a_1,\cdots ,a_{p-1},b}$ to the ``symbol'' 
$\{a_1,\cdots ,a_{p-1}\}\otimes (b)\bmod \mathcal{R}_p$, see \cite{ABa}
(or \cite{Baeza2} for a short exposition of these facts).

\begin{cor} \label{coho}
Let $x$ and $y$ be symbols in $H^p(k, \Zm2(p-1))$.  Suppose that there exist 
finite extensions $K_1, \ldots, K_n$ of $k$ such that $x$ and $y$ agree 
over every $K_i$.  Then there exists a Galois extension $K$ of $k$ such that 
$x$ and $y$ agree over $K$ and $[K:k]$ divides
   \[
    g:=\gcd\{4, [K_1:k], [K_2:k], \ldots, [K_n: k]\}.
      \]
If $p \ge 2$ and $g$ is even, the extension $K$ may be chosen to kill
both $x$ and $y$.  If $g=1$, then there exists a Galois extension of degree dividing $2$
that kills both $x$ and $y$. 
\end{cor}

\section{Proof of Theorem \ref{MT}: type $G_2$} \label{G2}

Every group $G$ of type $G_2$ and every $G$-torsor $X$
may be (in a compatible manner) identified with 
an octonion $k$-algebra \cite[33.24]{KMRT}.  Such algebras are determined by
their norm form, which is a 3-Pfister form \cite[\S1.7]{Sp:ex}. 

Suppose that the $G$-torsor $X$ is trivial (= has a closed point) over finite extensions $K_1, K_2, \ldots, K_n$ of $k$.  Let $\phi_G$ and $\phi_X$ be 3-Pfister forms corresponding to $G$ and $X$ respectively.  The quadratic forms are isomorphic over $K_i$ for all $i$ by hypothesis, and Prop.~\ref{quad.prop} gives a separable extension $K$ of $k$ such that $\phi_G$ and $\phi_X$ are $K$-isomorphic and $[K:k]$ divides $[K_i:k]$ for all $i$.  That is, $X$ has a $K$-point, and the answer to (\ref{ques}$^\prime$)---equivalently, \eqref{ques}---is ``yes'' for every group $G$ of type $G_2$ and every $G$-torsor $X$.

\section{Proof of Theorem \ref{MT}: type $F_4$} \label{F4}

An Albert $k$-algebra $J$ has $p$-Pfister form invariants $f_p(J)$ 
for $p = 3, 5$.  Every Albert algebra also has a ``mod 3'' invariant 
$g_3(J)$ living in the group $H^3(k, \Zm3(2))$.  (The correspondence between 3-Pfister
forms and  elements of $H^3(k, \Zm2(2))$ identifies the pair 
$(f_3(J),g_3(J))$ with an element of $H^3(k, \Zm6(2))$.  This element is the Rost
invariant of $J$ with respect to the split group of type $F_4$, up to sign.) 
The algebra $J$ is reduced if and only if $g_3(J)$ is zero \cite{Ro},
\cite{PR:el}.  There is a unique Albert algebra ---  the \emph{split} Albert
algebra --- with $f_3$ hyperbolic and $g_3$ equal to zero.

Let $G$ be a reduced group of $F_4$, i.e., $G$ is of the form $\Aut(J)$ 
for a reduced Albert $k$-algebra $J$.
We are given a $G$-torsor $X$ and finite extensions 
$K_1, K_2, \ldots, K_n$ of $k$ such that $X$ is trivial over 
$K_i$ for all $i$.  The torsor $X$ is the collection $\Iso(J', J)$ of 
isomorphisms $J' \iso  J$ for some Albert $k$-algebra $J'$, and
by hypothesis $J'$ and $J$ are $K_i$-isomorphic for all $i$.  We will 
construct a separable extension $K$ of $k$ such that $J'$ and $J$ are
$K$-isomorphic and the degree $[K:k]$ divides the gcd $g$ of the degrees
$[K_i : k]$. For $\ell = 2, 3$, we first construct an $\ell$-primary separable
extension $E_\ell$ of $k$ which makes the mod  $\ell$ invariants of $J$ and
$J'$ equal.

\smallskip

\underline{$\ell = 2$:} If $g$ is odd, then $[K_i : k]$ 
is odd for some $i$.  Hence $f_3$ and $f_5$ agree for $J$ and $J'$, 
and we take $E_2 = k$.  If $g$ is even, let $E_2$ be the Galois extension
provided by Prop.~\ref{quad.prop} that kills $f_3(J)$ and $f_3(J')$. 
Since $f_3$ divides $f_5$, the extension $E_2$ kills $f_5(J)$ and $f_5(J')$
also.

\smallskip

\underline{$\ell = 3$:} If $g_3(J')$ is zero, we may take $E_3 = k$.  
If $g$ is not divisible by 3, then $[K_i: k]$ is not divisible by 3 for some
$i$ and $g_3(J')$ is zero.  Otherwise, $J'$ is a division algebra and so contains 
a separable cubic subfield, which we take to be $E_3$.  Then $J' \ot_kE_3$ has 
zero divisors because $ E_3 \ot_k E_3$ does, whence $g_3(J')$ is killed by
$E_3$.

\smallskip

Set $K$ to be the compositum of $E_2$ and $E_3$ in some separable 
closure of $k$.  Since $E_2$ is Galois over $k$, the degree 
$[K: k]$ divides the product $[E_2 : k] [E_3 : k]$, hence $K$ has degree
dividing $g$. By construction, $K$ kills $g_3(J')$, so $J'$ is reduced over
$K$.  Since $f_i(J)$ agrees with $f_i(J')$ over $ K$ for $i = 3, 5$, the
algebras $J$ and $J'$ are isomorphic over $K$ by \cite[5.8.1]{Sp:ex},
\cite[4.1]{Ptr:struct}.


\section{The connecting homomorphism for $E_6$}

The purpose of this section is to prove a technical result we will need 
to prove Theorem \ref{MT} for groups of type $E_6$.

We say that a group $G$ \emph{has trivial Tits algebras} if the 
$k$-algebra $\End_G(V)$ is a field for every irreducible representation 
$V$ of $G$.  From here until the end of the paper, $H^1$ denotes flat 
cohomology as in \cite{DG} or \cite{Wa}.  
For smooth groups---in particular, for all the groups we consider here 
except for central subgroups of semisimple groups in bad 
characteristic---it agrees with the Galois cohomology
defined in \cite{SeCG}.

\begin{lem} \label{app.thm}
Let $C$ be a central subgroup of an algebraic group $G$ of type $E_6$.  
If $G$ is isotropic and has trivial Tits algebras, then the natural map
   \[
   \partial \!: (G/C)(k) \ra H^1(k, C)
   \]
is surjective.
\end{lem}

This lemma is an $E_6$-analogue of a well-known fact about classical groups.  
Specifically, let $G$ be the spin group of a nonsingular quadratic from $q$ on a 
vector space $V$.  Let $C$ be the kernel of the ``vector representation'' 
$G \ra GL(V)$; it is isomorphic to $\mmu2$.  The group $H^1(k, C)$ is 
identified with $\kxii$ by the Kummer exact sequence, and the connecting 
homomorphism $\partial$ is the spinor norm. When $q$ is isotropic, $\partial$ 
is surjective, as is well-known in quadratic form theory \cite[p.~78]{Baeza}.  
The proof of Lemma \ref{app.thm} below is easily adapted to give 
an alternative proof that the spinor norm is surjective.

The crux case of Lemma \ref{app.thm} is where $G$ is of type $\oEiso$, 
precisely the case we need in the rest of the paper.  When $k$ has good 
characteristic (i.e., not 2 or 3), the result in the crux case is an 
easy consequence of the theory of Albert algebras.  We give an 
algebraic-group-theoretic proof that is valid in all characteristics.

We will use without comment the fact that $\partial$ fits naturally into 
an exact sequence
 \[
 (G/C)(k) \xrightarrow{\partial} H^1(k, C) \xrightarrow{\iota} H^1(k, G).
 \]
In particular, $\partial$ is surjective if and only if $\iota$ is the zero map.

\begin{proof}[Proof of Lemma \ref{app.thm}]
Clearly we may assume that $C$ is not the trivial group, hence that $G$ is not adjoint, 
i.e., $G$ is simply connected.  Since the center of $G$ is a twisted form of $\mmu3$, 
the subgroup $C$ must be the entire center.

The case where $G$ is split is standard.  Indeed, $G$ contains a maximal 
$k$-torus $T$ that is $k$-split.  The center $C$ is contained in every maximal
torus, hence the map $\iota$ above factors through $H^1(k, T)$.  The group
$H^1(k, T)$ is zero by Hilbert 90, since $T$ is split.

\smallskip

Suppose now that $G$ is of type \oEiso, i.e., has Tits index \eqref{e6.ind}.  
Fix a maximal $k$-torus $T$ in $G$ containing a maximal $k$-split torus.  
Fix a set of simple roots $\D = \{ \alpha_1, \ldots, \alpha_6 \}$ of $G$ 
with respect to $T$, numbered as in \cite{Bou:g4}.  Let $\La$ (resp.\ 
$\La_r$) denote the weight (resp.\ root) lattice.  Let $\omega_i$ be the
fundamental weight corresponding to $\alpha_i$.  The absolute 
Galois group $\G$ of $k$ acts naturally on $\La$ and the fixed sublattice
$\La^\G$ consists of the weights $\omega$ such that  $(\alpha_i, \omega) = 0$
for all $i \ne 1, 6$ by \cite[p.~108, Cor.~6.9]{BoTi}.  In particular, the
weight $\omega_1$ lies in $\La^\G$.

The weight (resp.\ root) lattice is the group of cocharacters 
$\Gm \ra T/C$ (resp.\ $\Gm \ra T$), and the weight $\omega_1$ corresponds 
to a $k$-defined cocharacter $\lbar$.  Consulting the tables in the back of
\cite{ Bou:g4}, we find that $\omega_1$ is not in the root lattice, but
$3\omega_1$ is.  Let $\ell \!: \Gm \ra T$ be the $k$-cocharacter
corresponding to $3 \omega_1$.  We have a commutative diagram with exact rows 
 \[   \begin{CD}
  1 @>>> \mmu3 @>>> \Gm @>{3}>> \Gm @>>> 1 \\
  @. @VVV @V{\ell}VV @VV{\lbar}V @. \\
  1 @>>> C @>>> T @>>> T/C @>>> 1
  \end{CD}
  \]
Here, the arrow $\mmu3 \ra C$ arises because the composition 
$\mmu3 \ra \Gm \ra T/C$ is the zero map.
 Since $\omega_1$ is not in the root lattice, $\ell$ is an 
injection, ergo the map $\mmu3 \ra C$ is an injection.  
That is, the center $C$ of $G$ is contained in a rank 1 $k$-split torus,
namely the image of $\ell$.  The map $\iota$ is zero by the same argument as
in the case where $G$ is split. (This paragraph and the previous one are an
adaptation of the arguments behind the applications of the Gille-Merkurjev
Norm Principle in the last section of \cite{M:norm}.  The map $\mmu3 \ra C$
constructed here is $h(\lbar)$ in the notation of \cite[1.3]{M:norm}.)

\smallskip

If $G$ is of type $\oE$, the remaining possibility is that it has Tits index
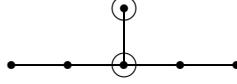
\[
\begin{picture}(5,2)
    \multiput(0.5,0.5)(1,0){5}{\circle*{\darkradE}}
    \put(2.5,1.5){\circle*{\darkradE}}

    \put(0.5,0.5){\line(1,0){4}}
    \put(2.5,1.5){\line(0,-1){1}}

    \put(2.5,0.5){\circle{\lradE}}
    \put(2.5,1.5){\circle{\lradE}}
\end{picture}
\]
We claim that this is impossible.  Indeed, the semisimple 
anisotropic kernel $G_\an$ of $G$ is of type
${^1\!A_2} \times {^1\!A_2}$.  Since the Tits algebras of $G$ are trivial, 
so are the Tits algebras of $G_\an$
by \cite[5.5.5]{Ti:R}, hence $G_\an$ is isotropic.  This is a contradiction.
 
\smallskip
 
Now suppose that $G$ is of type $\dE$ and let $K$ be the separable 
quadratic extension of $k$ over which $G$ is of type $\oE$.  The hypothesis on 
the Tits algebras ensures that $G$ is split or of type $\oEiso$ over $K$, 
hence that $C$ is contained in a $K$-defined rank 1 split torus $S$ in $G$.  
Not only is the map
   \[
     \partial_K \!: (G/C)(K) \ra H^1(K, C)
   \]
 surjective, it is even surjective when one restricts to the rational 
subgroup $S/C$ of $G/C$.
The Gille-Merkurjev Norm Principle from \cite{Gille:norm},
\cite{M:norm} shows that the image of $\partial$ contains the corestriction
$\cores_{K/k} H^1(K, C)$.  Since $K/k$ is quadratic and $C$ is 3-torsion, a
restriction-corestriction argument shows that the corestriction $H^1(K, C) \ra
H^1(k, C)$ is surjective. 
\end{proof}   

We remark that the Gille-Merkurjev
Norm Principle is proved for Galois cohomology and here we are using flat
cohomology, and the two cohomology theories may give different values for
$H^1(k, C)$ in characteristic 3.  However, the proof of the norm principle in
\cite{M:norm} goes through with no changes since flat and Galois cohomology
agree for reductive groups.  The only point requiring checking is Merkurjev's
Lemma 3.11, which is easy to translate.   

\begin{rmk}
The same sort of proof gives an $E_7$ version of Lemma \ref{app.thm}: 
Let $C$ be a central subgroup of   an algebraic group $G$ of type $E_7$.  
\emph{If $G$ is isotropic, has trivial Tits algebras, and is not of type 
$E^{66}_{7,1}$, then the natural map $(G/C)(k) \ra H^1(k, C)$ is surjective.}
\end{rmk}

\section{Proof of Theorem \ref{MT}: type $E_6$} \label{E6}

In this entire section, we assume that $G$ is a simply connected 
group of type $E_6$ that is split or of type \oEiso.  We write $Z$ 
for the center of $G$; it is isomorphic to the group scheme $\mmu3$ of cube
roots of unity.

\begin{lem} \label{small}
The group $G$ contains a subgroup $H$ that is reduced of type $F_4$, and the natural map
   \[
   H^1(k, H) \times H^1(k, Z) \ra H^1(k, G)
   \]
is surjective.
\end{lem}

\begin{proof}
Write $G_\sp$ for the split simply connected group of type $E_6$.  
The split group $H_\sp$ of type $F_4$ is contained in $G_\sp$, 
and the induced map
  \begin{equation} \label{incl}
  H^1(k, H_\sp) \times H^1(k, Z) \ra H^1(k, G_\sp)
  \end{equation}
is surjective.  This has been shown in \cite[3.4]{G:rinv} under the assumption that 
$k$ has characteristic different from 2, 3,
but this restriction is unnecessary by the following reasoning.  
In all characteristics $G_\sp $
is the group of isometries of a cubic form---a norm for an Albert algebra
$J$---and the elements of norm 1 in  $J$ form an open orbit in $\P(J)$
\cite[3.16(3)]{Asch:E6}.   The stabilizer of the identity element in $J$ is
of type $F_4$ in all characteristics by \cite[4.6]{Sp:jord}, hence the
stabilizer of the identity element in $\P(J)$ is $F_4 \times \mmu3$.  The
conclusion now follows by \cite[Lemma 3.1]{G:rinv}, where the last paragraph
of the proof is replaced with an appeal to \cite[p.~373, Prop.~III.4.4.6b]{DG},
a flat cohomology analogue of \cite[I.5.4, Prop.~37]{SeCG}.  
This completes the proof when $G$ is split.

Now suppose that $G$ is of type $\oEiso$.  The semisimple anisotropic 
kernel $G_\an$ is simply connected of type $D_4$.  It has trivial 
Tits algebras \cite[5.5.5]{Ti:R}, 
so is obtained by twisting the split simply connected group 
$\Spin_8$ of type $D_4$ by a 1-cocycle $\alpha$ with values in $\Spin_8$. 
The group $\Spin_8$ is even a subgroup of $H_\sp$ in $G_\sp$.  Twisting the
split group $G_\sp$ by $\alpha$ produces a group of type $\oE$ with the same
semisimple anisotropic kernel as $G$, hence the group is isomorphic to $G$ by
Tits' Witt-type theorem.  Twisting $H_\sp$ by $\alpha$, we find a group $H$
of type $F_4$ that is reduced.  Twisting everything in  \eqref{incl} by
$\alpha$ gives the lemma. 
\end{proof}

Let $r_G$ denote the Rost invariant of $G$ as defined in \cite{MG}.  The composition
\begin{equation} \label{rH}
\begin{CD}
H^1(k, H) @>>> H^1(k, G) @>{r_G}>> H^3(k, \QZt),
\end{CD}
\end{equation}
is equal to an integer multiple $n r_H$ of the Rost invariant 
$r_H \!: H^1(k, H) \ra H^3(k, \QZt)$.  Over a separable closure of $k$, 
the subgroup $H$ from Lemma \ref{small} is similar to the standard inclusion 
of the split $F_4$ into the split $E_6$, and the integer multiple for this 
inclusion is 1 by \cite[2.4]{G:rinv}.  Since this multiplying factor does not 
change under scalar extensions \cite[Prop.~7.9(4)]{MG}, the 
composition \eqref{rH} is equal to $r_H$.

By Lemma \ref{app.thm}, the image of the map $H^1(k, Z) \ra H^1(k, G)$ is zero.  
Since the Rost invariant is compatible with twisting, for 
$(\alpha, \la) \in H^1(k, H) \times H^1(k, Z)$ we have:
\begin{equation} \label{rGrH}
r_G(\alpha, \la) = r_H(\alpha).
\end{equation}

\begin{prop} \label{rinv}
The Rost invariant $H^1(k, G) \ra H^3(k, \QZt)$ has kernel zero.
\end{prop}

\begin{proof}
Suppose $\beta$ is in the kernel.
It is the image of some pair $(\alpha, \la)$ in $H^1(k, H) \times H^1(k, Z)$ 
such that $r_H(\alpha)$ is zero by \ref{small} and \eqref{rGrH}.  

The group $H$ is the automorphism group of a reduced Albert algebra $J$.  
Let $J'$ be the Albert algebra corresponding to $\alpha$.  
Since $r_H(\alpha)$ is zero, the Rost invariants of $J$ and $J'$ relative 
to the split group of type $F_4$ agree, i.e., $f_3(J) = f_3(J')$ and 
$g_3(J) = g_3(J')$.  In particular, $J'$ is also reduced.  
An explicit formula for the cubic norm form on a reduced Albert algebra 
is well-known, and since the 3-Pfister forms $f_3(J)$ and $f_3(J')$ are the same, 
it is easy to cook up an isometry between the norms on $J$ and $J'$.  
Since $G$ may be viewed as the group of isometries of the cubic norm on $J$, 
the image of $\alpha$ in $H^1(k, G)$ is zero by descent.

The group $H^1(k, Z)$ acts on $H^1(k, G)$, and $\beta$ is 
$\la \cdot (\im \alpha)$.  Hence $\beta$ is in the image of the map 
$H^1(k, Z) \ra H^1(k, G)$.  But this map has image zero by Lemma \ref{app.thm},
hence $\beta$ is zero.
\end{proof}

Now we may prove Theorem \ref{MT} for our $G$.  
By hypothesis, we are given a $G$-torsor $X$ and finite extensions 
$K_1, K_2, \ldots, K_n$ of $k$ such that $X$ is trivial over $K_i$ for all
$i$.  Fix a pair $(\alpha, \la)$  in $H^1(k, H) \times H^1(k, Z)$ mapping to
$X$.  Then $r_H(\alpha)$ equals $r_G(X)$, and it is killed by $K_i$ for all $i$.  
Let $K$ be the
extension constructed from $H$, $\alpha$, and the $K_i$ as in \S\ref{F4}.  The
extension $K$ kills $r_G(X)$, hence trivializes $X$ by the proposition.  By
construction, $K$ is separable over $k$ of degree dividing $[K_i : k]$ for all
$i$.  This completes  the proof of the theorem.

\section{Final remarks} \label{specsec}

In the special case where $d = 1$ and $X$ is a $G$-torsor, 
Totaro's question overlaps an earlier question by 
Serre \cite[p.~233]{SeCG:p}.  But even this simpler case is wide open for
general groups of type $F_4$. 
\emph{For $k$ of characteristic $\ne 2, 3$, the following are equivalent:}
\begin{align}
&\parbox{4in}{\emph{For every $k$-group $G$ of type $F_4$, every $G$-torsor $X$, 
and $d = 1$, 
the answer to \eqref{ques} is ``yes''.}} \label{rost.1} \\
&\parbox{4in}{\emph{Albert $k$-algebras are classified by their invariants 
$f_3$, $f_5$, and $g_3$.}} \label{rost.2}
\end{align}
Indeed, the equivalence of \eqref{rost.1} and \eqref{rost.2} is an easy 
consequence of a theorem of Rost from \cite{Rost:albert}, which says: 
If $J$, $J'$ are Albert $k$-algebras such that $f_3$, $f_5$, and $g_3$ 
agree for $J$ and $J'$, then there exist extensions $K_1$, $K_2$ such 
that $J$, $J'$ are isomorphic over $K_1$ and $K_2$ and such that the degrees 
$[K_1 : k]$ and $[K_2 : k]$  are coprime.
But proving or disproving \eqref{rost.2} is Question \#1 on Petersson-Racine's
list of open problems in \cite{PR}; it is viewed by many as the principal 
outstanding problem concerning Albert algebras.  

\smallskip

Suppose for the moment that \eqref{rost.1} is true.  
To extend \eqref{rost.1} to the case of general $d$, 
one would want an analogue of Corollary \ref{coho} for the prime 3.

\smallskip

If the answer to Totaro's question is yes for 
$G$ simply connected of type $\oE$ with trivial Tits algebras, then one can 
deduce (in a manner very similar to the above)
that two Albert algebras are isotopic if and only if they have the same $f_3$
and $g_3$ invariants.  Whether or not that is true is Question \#4 in
Petersson-Racine's list.

{\small 
\subsection*{Acknowledgements} The
first author thanks FernUniversit\"at Hagen for its hospitality while a first draft of this article was written, and Holger Petersson for the particularly simple argument for the prime 3 in \S\ref{F4}.

The second author has been supported in part by the European research network
``Algebraic $K$-Theory, Linear Algebraic Groups and Related Structures''
HPRN-CT-2002-00287.  He thanks the University of Ghent for its hospitality during
a stay where some of the work on this article has been carried out.}

\providecommand{\bysame}{\leavevmode\hbox to3em{\hrulefill}\thinspace}
\providecommand{\MR}{\relax\ifhmode\unskip\space\fi MR }
\providecommand{\MRhref}[2]{%
  \href{http://www.ams.org/mathscinet-getitem?mr=#1}{#2}
}
\providecommand{\href}[2]{#2}

\end{document}